\documentclass[11pt,twoside]{article}
\usepackage{amsmath,latexsym,amssymb,amsfonts,amsbsy}
\usepackage[mathscr]{eucal} 
\usepackage{graphicx}

\newfont{\bb}{msbm10}

\newtheorem{definition}{Definition}[section]

\newtheorem{proposition}[definition]{Proposition}

\newtheorem{corollary}[definition]{Corollary}

\begin{document}
\cleardoublepage \pagestyle{myheadings}

\bibliographystyle{plain}

\title{\bf Tensor Transpose and Its Properties}

\author{Ran Pan
\\
\
{\it Department of Mathematics}\\
{\it University of California, San Diego}\\
{\it Email: r1pan@ucsd.edu }\\[2mm]
}

\maketitle \markboth{\small  Ran Pan}
{\small  Tensor Transpose and Its Properties}

\begin{abstract}
Tensor transpose is a higher order generalization of matrix transpose. In this paper, we use permutations and symmetry group to define the tensor transpose. Then we discuss the classification and composition of tensor transposes. Properties of tensor transpose are studied in relation to tensor multiplication, tensor eigenvalues, tensor decompositions and tensor rank.
\end{abstract}

\bigskip


\bigskip

\bigskip

\section{Introduction} \label{intro-sec}
A tensor is a multidimensional or \emph{N}-way array. The order of a tensor is the number of its dimensions. If $\boldsymbol{\mathscr{X}}$ denotes a real \emph{N}-order tensor, we have $\boldsymbol{\mathscr{X}}\in \mathbb{R}^{ I_{1}\times{I_{2}}\times{\cdots}\times{I_{N}}} $. $I_{1}\times{I_{2}}\times{\cdots}\times{I_{N}}$ is called the size of $\boldsymbol{\mathscr{X}}$. It's clear that tensor $\boldsymbol{\mathscr{X}}$ with entries $\boldsymbol{\mathscr{X}}(i_{1},i_{2},\ldots,i_{N})$ has \emph{N} indices. There're some examples: a vector (i.e., a 1-order tensor), a matrix (i.e., a 2-order tensor)and a 3-order tensor are shown in Figure 1.1.
\begin{figure}[!ht]
\centering
\includegraphics[bb= 0 0 400 200, width=100mm, height=50mm]{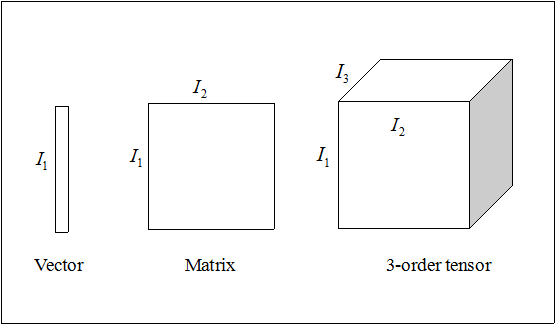}
\caption{1-order, 2-order and 3-order tensors}
\end{figure}
\par
In recent decades, research of tensors attracted much attention. In theoretical work, the theories of tensor multiplication and decompositions are much developed, as well as eigenvalues and singular values of tensors~\cite{Qi2}~\cite{Qi3}~\cite{LHLim}. For applications, tensors appear in many fields, such as psychology~\cite{psy}, web searching~\cite{TOPHITS} and so on. Although the notion of tensor transpose is often mentioned together with supersymmetric tensors ~\cite{STSTR}, specific discussions concerning tensor transpose draw less attention.\par
In this paper, tensors are viewed from a different prospective: tensor transpose. It is known that matrix transpose is a notion of a matrix and has several properties. Here, we focus our emphases on extending the notion of transpose from 2-order tensors to high order tensors and discovering its properties. \par
In Section 2, we will review knowledge of permutations and symmetry group~\cite{Jacobson}, and use them to define the tensor transpose. Then we discuss the classification and composition of tensor transposes. In Section 3, the relationship between transpose and tensor multiplication will be discussed. It is known transpose of matrix multiplication satisfies that $A^{T}B^{T}=(BA)^{T}$, where $A$ and $B$ are matrices. Proposition 3.1, as the most important result of this paper, will be introduced. In section 4, we prove that tensor $l^{p}$-eigenvalues are invariant under some certain tensor transpose. In Section 5, the discussion about the transpose and decomposition will be continued. Two methods of decompositions CP decomposition and Tucker decomposition are considered separately in corresponding passages respectively. \par

\section{Definition of Tensor Transpose} \label{Not-sec}
 Before introducing the definition of tensor transpose, we review fundamental knowledge of permutations and symmetry group first.
\par

Assume set $S=\{1,2,3,\cdots,n\}$ and $\sigma$ is a permutation of $S$, usually, $\sigma$ is denoted as follows,  with $\sigma(i)=k_{i}$
\begin{equation*}
  \sigma=\left(
           \begin{array}{ccccc}
             1 & 2 & 3 & \cdots & n \\
             k_{1} & k_{2} & k_{3} & \cdots & k_{n} \\
           \end{array}
         \right),
\end{equation*}where $\{k_{1},k_{2},k_{3},\cdots,k_{n}\}=\{1,2,3,\cdots,n\}=S$.The same as functions, two permutations can be composed. For example, $\sigma$ and $\tau$ are two permutations and their composition $\tau\sigma(i)=\tau(\sigma(i))$.

If there exist $r$ different numbers $i_{1}$,$i_{2}$,$\cdots$,$i_{r}$ subjected to
\begin{equation*} \sigma(i_{1})=i_{2},\sigma(i_{2})=i_{3},\cdots,\sigma(i_{r-1})=i_{r},\sigma(i_{r})=i_{1},\end{equation*}
the permutation $\sigma$ is called a r-circle.

All permutations are products of cycles. Specially, a $3$-order permutation can be written in form of a cycle. For instance,
\begin{equation*}
  \sigma=\left(
           \begin{array}{ccc}
             1 & 2 & 3  \\
             2 & 1 & 3  \\
           \end{array}
         \right)
         =\left(\begin{array}{cc}
            1 & 2
          \end{array}\right)
         .
\end{equation*}
 We denote $\sigma^{-1}$ as the inverse of $\sigma$.
Then we review the definition of a symmetry group.
$S$ is a finite set consisting of $n$ elements. The set consisting of all permutations of any given set $S$, together with the composition of function is symmetry group of $S$. The notation of symmetry group is $S_{n}$.
 Apparently, there are $n!$ permutations for the set $S$. 

 Now we can introduce the definition of tensor transpose. We know that matrix transpose is a permutation of the two indices. Considering the fact that a matrix is a 2-order tensor, in this paper, we extend the definition of transpose to high order tensors.
\begin{definition}
Let $\boldsymbol{\mathscr{X}}$ be an $n$-order tensor. $\boldsymbol{\mathscr{Y}}$ is called tensor transpose of $\boldsymbol{\mathscr{X}}$ associated with $\sigma$, if entries $\boldsymbol{\mathscr{Y}}(i_{\sigma(1)},i_{\sigma(2)},\ldots,i_{\sigma(n)})=\boldsymbol{\mathscr{X}}(i_{1},i_{2},\ldots,i_{n})$, where $\sigma$ is an element of $S_{n}$ but not an identity permutation. $\boldsymbol{\mathscr{Y}}$ is denoted by $\boldsymbol{\mathscr{X}}^{T\sigma}$.
\end{definition}
This definition shows if $\boldsymbol{\mathscr{X}}$ is of size $I_{1}\times{I_{2}}\times{\cdots}\times{I_{N}}$, $\boldsymbol{\mathscr{X}}^{T\sigma}$ is
 of size $I_{\sigma(1)}\times{I_{\sigma(2)}}\times{\cdots}\times{I_{\sigma(N)}}$. Elementwise, for example, assume $\boldsymbol{\mathscr{X}}$ is a 3-order tensor and $\sigma=\left(\begin{array}{ccc}
                                                                          1 & 2 & 3
                                                                        \end{array}\right)\in S_{3}$, i.e. $\sigma(1)=2, \sigma(2)=3, \sigma(3)=1
$. We have $\boldsymbol{\mathscr{X}}^{T\sigma}(i_{2},i_{3},i_{1})=\boldsymbol{\mathscr{X}}(i_{1},i_{2},i_{3})$.\par
Using Bader and Kolda's MATLAB tensor toolbox~\cite{toolbox}, a tensor can be transposed in MATLAB. In tensor toolbox, the function PERMUTE is used to transpose a tensor.\\
\par

In order to facilitate the description and discussion, some special symbols are utilized for transpose of 3-order tensors.
\begin{equation*}
\begin{cases}
\text{positive transpose}\phantom{12}\boldsymbol{\mathscr{X}}^{T_{+}}=\boldsymbol{\mathscr{X}}^{T\sigma_{1}} & \text{where $\sigma_{1}=\left(\begin{array}{ccc}
            1 & 2 &3
          \end{array}\right)$}\\
         \text{negative transpose }\phantom{1}\boldsymbol{\mathscr{X}}^{T_{-}}=\boldsymbol{\mathscr{X}}^{T\sigma_{2}} & \text{where $\sigma_{2}=\left(\begin{array}{ccc}
            1 & 3 &2
          \end{array}\right)$}\\
         \text{first transpose }\phantom{1211}\boldsymbol{\mathscr{X}}^{T_{1}}=\boldsymbol{\mathscr{X}}^{T\sigma_{3}} & \text{where $\sigma_{3}=\left(\begin{array}{cc}
            2 &3
          \end{array}\right)$}\\
        \text{second transpose}\phantom{11}\boldsymbol{\mathscr{X}}^{T_{2}}=\boldsymbol{\mathscr{X}}^{T\sigma_{4}} & \text{where $\sigma_{4}=\left(\begin{array}{cc}
            1  &3
          \end{array}\right)$}\\
        \text{third transpose } \phantom{112}\boldsymbol{\mathscr{X}}^{T_{3}}=\boldsymbol{\mathscr{X}}^{T\sigma_{5}} & \text{where $\sigma_{5}=\left(\begin{array}{cc}
            1 & 2
          \end{array}\right)$}\\
\end{cases}
\end{equation*}
Considering a 3-order tensor is like a rectangular cuboid, transpose of a 3-order tensor has its geometric meaning. Transpose can be looked upon as a rotation of the rectangular cuboid.\par
In terms of an $n$-order tensor, it is clear that a certain transpose is corresponding to a certain permutation. Therefore, we've got the first property of tensor transpose in the paper.
\begin{proposition}
An $n$-order tensor has $n!-1$ different transposes.
\end{proposition}\par
In addition, according to definition of transpose, we can define the supersymmetric tensor in a new way. Symmetry is an important notion for matrices, while it is called supersymmetry for high order tensors. In scientists' previous work, a supersymmetric tensor is often described as a tensor whose entries are invariant under any permutation of their indices. We know that matrix $A$ is symmetric if $A^{T}=A$. A similar definition of supersymmetric tensor is given as follows.
\begin{definition}
$\boldsymbol{\mathscr{X}}$ is called a supersymmetric tensor, if $\boldsymbol{\mathscr{X}}=\boldsymbol{\mathscr{X}}^{T\sigma}$, for all $\sigma\in{S_{n}}$, where $n$ is the order of $\boldsymbol{\mathscr{X}}$.
\end{definition}

From the viewpoint of group theory, tensor transpose can be treated as an action of a symmetry group on a tensor and supersymmetric tensors can be treated as fixed elements of group $S_{n}$.

Besides, according to the definition of tensor transpose, tensor transposes can be classified into two classes: total transpose and partial transpose.
 \begin{definition}$\boldsymbol{\mathscr{X}}^{T\sigma}$ is called a total transpose, if $\sigma$ is a derangement. Otherwise, it's a partial transpose. A derangement~\cite{DERA} is a permutation such that none of the elements appear in their original position. It means $\sigma$ is derangement if  $\sigma(i)\neq{i}$, for all ${i}\in\{1,2,\cdots,n\}$.
  \end{definition}
  Assume an $n$-order tensor has $a_{n}$ total transposes, apparently $a_{n}$ is also the number of derangements in symmetry group $S_{n}$. We have
\begin{equation*}
a_{n}=n!\sum^{n}_{i=0}{\frac{(-1)^{i}}{i!}},
\end{equation*}
and it is called "de Montmort number"~\cite{deMon}.
  For 3-order tensor $\boldsymbol{\mathscr{X}}$, $a_{3}=2$, it has two total transposes $\boldsymbol{\mathscr{X}}^{T_{+}}$ and $\boldsymbol{\mathscr{X}}^{T_{-}}$. And $\boldsymbol{\mathscr{X}}^{T_{1}}$, $\boldsymbol{\mathscr{X}}^{T_{2}}$ and $\boldsymbol{\mathscr{X}}^{T_{3}}$ are partial transpose.\par

It is known that two permutations can be composed, as well as two tensor transposes. For matrix $A$, $(A^{T})^{T}=A$. For high order tensors, we have following proposition.
\begin{proposition}
$(\boldsymbol{\mathscr{X}}^{T\sigma})^{T\tau}=\boldsymbol{\mathscr{X}}^{T(\tau\sigma)}$, where $\tau\sigma$ means the composition of the two permutations.
\end{proposition}

From this proposition, it is obvious that the composition of transpose is equivalent to composition of permutations. If $\tau=\sigma^{-1}$, $(\boldsymbol{\mathscr{X}}^{T\sigma})^{T\sigma^{-1}}=\boldsymbol{\mathscr{X}}$. Therefore, we have the property for 3-order tensors.
\begin{corollary}Let $\boldsymbol{\mathscr{X}}$ be a 3-order tensor.
$(\boldsymbol{\mathscr{X}}^{T_{+}})^{T_{-}}=\boldsymbol{\mathscr{X}}$, $(\boldsymbol{\mathscr{X}}^{T_{-}})^{T_{+}}=\boldsymbol{\mathscr{X}}$, $(\boldsymbol{\mathscr{X}}^{T_{1}})^{T_{1}}=\boldsymbol{\mathscr{X}}$, $(\boldsymbol{\mathscr{X}}^{T_{2}})^{T_{2}}=\boldsymbol{\mathscr{X}}$, $(\boldsymbol{\mathscr{X}}^{T_{3}})^{T_{3}}=\boldsymbol{\mathscr{X}}$.
\end{corollary}\par

\section{Tensor Transpose of Tensor Multiplication} \label{TTM-sec}
In this section, we will consider tensor transpose of tensor multiplication. High order tensor multiplication is much more complex than matrix multiplication, and high order tensors have more transposes than matrices (i.e. 2-order tensors). A full treatment of tensor multiplication can be found in Bader and Kolda's work ~\cite{MOFHOD}~\cite{MTC}. Here we only discuss some kinds multiplication of them.\par

\subsection{Tensor matrix multiplication}
A familiar property of matrix transpose and matrix multiplication is that $A^{T}B^{T}=(BA)^{T}$, where $A$ and $B$ are matrices. However, a high order tensor usually has more than three dimensions. Therefore, we must specify which dimension is multiplied by the matrix. In this passage, we adopt $n$-mode product ~\cite{LDL}.\par
Let $\boldsymbol{\mathscr{X}}$ be an $I_{1}\times{I_{2}}\times{\cdots}\times{I_{N}}$ tensor and $\textbf{U}$ be a $J_{n}\times{I_{n}}$ matrix. Then the $n$-mode product of $\boldsymbol{\mathscr{X}}$ and $\textbf{U}$ is denoted by ${\boldsymbol{\mathscr{X}}}\times_{n}\textbf{U}$ and its result is a tensor of size ${{I_{1}}\times{I_{2}}{\times}{\cdots}{\times}{I_{n-1}}\times{J_{n}}{\times}{I_{n+1}}{\times}{\cdots}\times{I_{N}}}$. The element of ${\boldsymbol{\mathscr{X}}}\times_{n}\textbf{U}$ is defined as
\begin{equation*}
({\boldsymbol{\mathscr{X}}}\times_{n}\textbf{U})({i_{1}},{\ldots},{i_{n-1}},{j_{n},{i_{n+1}},{\ldots},{i_{N}}})=\sum_{i_{n}=1}^{I_{n}}{\boldsymbol{\mathscr{X}}({i_{1}},{i_{2}},{\ldots},i_{N})}\textbf{U}({j_{n}},{i_{n}}).
\end{equation*}
\begin{proposition}
         Let ${\boldsymbol{\mathscr{X}}}$ be an $N$-order tensor, we have \begin{equation*}{\boldsymbol{\mathscr{X}}}^{T\sigma}{\times_{\sigma^{-1}(n)}}{\rm{\bf{U}}}=({\boldsymbol{\mathscr{X}}}{\times_{n}}{\rm\bf{U}})^{T\sigma}.\end{equation*}
\end{proposition}
Proof:\\
Let ${\boldsymbol{\mathscr{Y}}}={\boldsymbol{\mathscr{X}}}^{T\sigma}$, ${\boldsymbol{\mathscr{Z}}}={\boldsymbol{\mathscr{Y}}}{\times_{\sigma^{-1}(n)}}{\rm{\bf{U}}}$, ${\boldsymbol{\mathscr{P}}}={\boldsymbol{\mathscr{X}}}{\times_{(n)}}{\rm{\bf{U}}}$, ${\boldsymbol{\mathscr{Q}}}={\boldsymbol{\mathscr{P}}}^{T\sigma}$.\\
According to the definition of tensor transpose, we have
\begin{equation*}
{\boldsymbol{\mathscr{Y}}}(i_{1},i_{2},\ldots,i_{N})={\boldsymbol{\mathscr{X}}}(i_{\sigma^{-1}(1)},i_{\sigma^{-1}(2)},\ldots,i_{\sigma^{-1}(N)}),
\end{equation*}
and
\begin{equation*}
{\boldsymbol{\mathscr{Z}}}({i_{1}},{\ldots},{i_{n-1}},{j_{n},{i_{n+1}},{\ldots},{i_{N}}})={\sum_{i_{\sigma^{-1}(n)}=1}^{I_{{\sigma^{-1}(n)}}}}{\boldsymbol{\mathscr{Y}}({i_{1}},{i_{2}},{\ldots},i_{N})}\textbf{U}({j_{\sigma^{-1}(n)}},{i_{\sigma^{-1}(n)}}),
\end{equation*}
that is
\begin{equation*}
{\boldsymbol{\mathscr{Z}}}({i_{1}},{\ldots},{i_{n-1}},{j_{n},{i_{n+1}},{\ldots},{i_{N}}})={\sum_{i_{\sigma^{-1}(n)}=1}^{I_{{\sigma^{-1}(n)}}}}{\boldsymbol{\mathscr{X}}}(i_{\sigma^{-1}(1)},i_{\sigma^{-1}(2)},\ldots,i_{\sigma^{-1}(N)})\textbf{U}({j_{\sigma^{-1}(n)}},{i_{\sigma^{-1}(n)}}).
\end{equation*}
For the right side of the equation,
\begin{equation*}
{\boldsymbol{\mathscr{P}}}({i_{1}},{\ldots},{i_{n-1}},{j_{n},{i_{n+1}},{\ldots},{i_{N}}})=\sum_{i_{n}=1}^{I_{n}}{\boldsymbol{\mathscr{X}}({i_{1}},{i_{2}},{\ldots},i_{N})}\textbf{U}({j_{n}},{i_{n}}),
\end{equation*}
since ${\boldsymbol{\mathscr{Q}}}={\boldsymbol{\mathscr{P}}}^{T\sigma}$,
\begin{equation*}
{\boldsymbol{\mathscr{Q}}}(i_{1},i_{2},\ldots,i_{N})={\boldsymbol{\mathscr{P}}}(i_{\sigma^{-1}(1)},i_{\sigma^{-1}(2)},\ldots,i_{\sigma^{-1}(N)}),
\end{equation*}
\begin{equation*}
{\boldsymbol{\mathscr{Q}}}({i_{1}},{\ldots},{i_{n-1}},{j_{n},{i_{n+1}},{\ldots},{i_{N}}})={\sum_{i_{\sigma^{-1}(n)}=1}^{I_{{\sigma^{-1}(n)}}}}{\boldsymbol{\mathscr{X}}}(i_{\sigma^{-1}(1)},i_{\sigma^{-1}(2)},\ldots,i_{\sigma^{-1}(N)})\textbf{U}({j_{\sigma^{-1}(n)}},{i_{\sigma^{-1}(n)}}).
\end{equation*}
Therefore, ${\boldsymbol{\mathscr{Z}}}({i_{1}},{\ldots},{i_{n-1}},{j_{n},{i_{n+1}},{\ldots},{i_{N}}})={\boldsymbol{\mathscr{Q}}}({i_{1}},{\ldots},{i_{n-1}},{j_{n},{i_{n+1}},{\ldots},{i_{N}}})$. $\blacksquare$
\\
\par

In fact, Proposition 3.1 is an extension of 2-order situation: $A^{T}B^{T}=(BA)^{T}$. For matrices, $A{\times}_{1}B=BA$, $A{\times}_{2}B=AB^{T}$. According to proposition 3.1, $(BA)^{T}=(A{\times}_{1}B)^{T}={A^{T}}{\times}_{2}B=A^{T}B^{T}$.
\subsection{Tensor inner product}
We consider three general scenarios for tensor-tensor multiplication: outer product,
contracted product, and inner product~\cite{MTC}.
\\
For outer product and contracted product, tensor transposes of them do not have significant characteristics.\\
For the inner product of two tensors, it requires that these two tensors are of the same size. Assume $\boldsymbol{\mathscr{A}}$, $\boldsymbol{\mathscr{B}}$ are two tensors of size $I_{1}{\times}I_{2}{\times}\cdots{\times}I_{N}$, the inner product of $\boldsymbol{\mathscr{A}}$, $\boldsymbol{\mathscr{B}}$ is given by
\begin{equation*}
\langle \boldsymbol{\mathscr{A}}, \boldsymbol{\mathscr{B}}\rangle={\sum_{i_{1}=1}^{I_{1}}}{\sum_{i_{2}=1}^{I_{2}}}\cdots{\sum_{i_{N}=1}^{I_{N}}}\boldsymbol{\mathscr{A}}(i_{1},i_{2},\ldots,i_{N})\boldsymbol{\mathscr{B}}(i_{1},i_{2},\ldots,i_{N})
\end{equation*}
\begin{proposition}
$\boldsymbol{\mathscr{A}}$ and $\boldsymbol{\mathscr{B}}$ are two tensors of the same size, we have $\langle \boldsymbol{\mathscr{A}}^{T\sigma}, \boldsymbol{\mathscr{B}}^{T\sigma}\rangle=\langle \boldsymbol{\mathscr{A}}, \boldsymbol{\mathscr{B}}\rangle.$
\end{proposition}
This property follows the definition of the inner product directly.
\par
Using the inner product, the Frobenius norm of a tensor is given by $\|\boldsymbol{\mathscr{A}}\|_{F}=\sqrt{\langle \boldsymbol{\mathscr{A}}, \boldsymbol{\mathscr{A}}\rangle}$. Then we have the following property.
\begin{corollary}
$\|\boldsymbol{\mathscr{A}}^{T\sigma}\|_{F}=\|\boldsymbol{\mathscr{A}}\|_{F}$
\end{corollary}

\section{Eigenvalues of Transposed Tensors}
It is known that eigenvalues keep invariant after the matrix is transposed. There are similar situations for eigenvalues of transposed tensors. In this section, we adopt Lim's definition of $l^{p}$-eigenvalues of nonsymmetric tensors ~\cite{LHLim}.
\\
\par
$\boldsymbol{\mathscr{A}}\in \mathbb{R}^{ n\times{n}\times{\cdots}\times{n}}$ is a $K$-order tensor. The homogeneous polynomial associated with tensor $\boldsymbol{\mathscr{A}}$ can be conveniently expressed as
\begin{equation*}
\boldsymbol{\mathscr{A}}(\textbf{x},\cdots,\textbf{x}):=\boldsymbol{\mathscr{A}}\times_{1}\textbf{x}\cdots\times_{K}\textbf{x}.
\end{equation*}
Since $\boldsymbol{\mathscr{A}}$ has $K$ sides, the tensor has $K$ different forms of eigenpairs as follows
\begin{equation*}
\boldsymbol{\mathscr{A}}(I_{n},\textbf{x}_{1},\cdots,\textbf{x}_{1})=\lambda_{1}\varphi_{p-1}(\textbf{x}_{1})
\end{equation*}
\begin{equation*}
\boldsymbol{\mathscr{A}}(\textbf{x}_{2},I_{n},\cdots,\textbf{x}_{2})=\lambda_{2}\varphi_{p-1}(\textbf{x}_{2})
\end{equation*}
\begin{equation*}
\cdots
\end{equation*}
\begin{equation*}
\boldsymbol{\mathscr{A}}(\textbf{x}_{K},\textbf{x}_{K},\cdots,I_{n})=\lambda_{K}\varphi_{p-1}(\textbf{x}_{K})
\end{equation*}
where $I_{n}$ is an $n$-by-$n$ identity matrix, $\varphi_{p-1}(\textbf{x}):=[sgn(x_{1}){|x_{1}|}^{p},\cdots,sgn(x_{n}){|x_{n}|}^{p}]$, and $sgn(x)$ is the sign function. The unit vector $\textbf{x}_{i}$ is called mode-$i$ eigenvector of $\boldsymbol{\mathscr{A}}$ corresponding to the mode-$i$ eigenvalue $\lambda_{i}$, $i=1,2,\cdots,K$.
\begin{proposition}
The mode-{i} eigenpairs are invariant under tensor transopose associated with $\sigma$, if $\sigma(i)=i$.
\end{proposition}
Proof:\\
Here, we take mode-$1$ as an example.\\
$\boldsymbol{\mathscr{A}}\in \mathbb{R}^{ n\times{n}\times{\cdots}\times{n}}$ is a $K$-order tensor. Let $\boldsymbol{\mathscr{B}}=\boldsymbol{\mathscr{A}}^{T{\sigma}}$, where $\sigma(1)=1$.\\
We have
\begin{equation*}\boldsymbol{\mathscr{B}}(I_{n},\textbf{x},\cdots,\textbf{x})=\boldsymbol{\mathscr{B}}\times_{2}\textbf{x}\times_{3}\textbf{x}\cdots\times_{K}\textbf{x}
=\boldsymbol{\mathscr{A}}^{T\sigma}\times_{2}\textbf{x}\times_{3}\textbf{x}\cdots\times_{K}\textbf{x}\end{equation*}
According to Proposition 3.1  ${\boldsymbol{\mathscr{X}}}^{T\sigma}{\times_{n}}{\rm{\bf{U}}}=({\boldsymbol{\mathscr{X}}}{\times_{\sigma(n)}}{\rm\bf{U}})^{T\sigma}$ and $\sigma(1)=1$,
\begin{equation*}
\begin{split}
\boldsymbol{\mathscr{A}}^{T\sigma}\times_{2}\textbf{x}\times_{3}\textbf{x}\cdots\times_{K}\textbf{x}
& =(\boldsymbol{\mathscr{A}}\times_{\sigma(2)}\textbf{x})^{T\sigma}\times_{3}\textbf{x}\cdots\times_{K}\textbf{x}\\
& =(\boldsymbol{\mathscr{A}}\times_{\sigma(2)}\textbf{x}\times_{\sigma(3)}\textbf{x})^{T\sigma}\cdots\times_{K}\textbf{x}\\
& \cdots\\
& =(\boldsymbol{\mathscr{A}}\times_{\sigma(2)}\textbf{x}\times_{\sigma(3)}\textbf{x}\cdots\times_{\sigma(K)}\textbf{x})^{T\sigma}\\
& = (\boldsymbol{\mathscr{A}}\times_{2}\textbf{x}\times_{3}\textbf{x}\cdots\times_{K}\textbf{x})^{T\sigma}\\
& =\boldsymbol{\mathscr{A}}\times_{2}\textbf{x}\times_{3}\textbf{x}\cdots\times_{K}\textbf{x}
\end{split}
\end{equation*}
Therefore, $\boldsymbol{\mathscr{B}}(I_{n},\textbf{x},\cdots,\textbf{x})=\boldsymbol{\mathscr{A}}(I_{n},\textbf{x},\cdots,\textbf{x})$.
When $\sigma(1)=1$, $\boldsymbol{\mathscr{A}}^{T\sigma}$ has the same eigenpairs with $\boldsymbol{\mathscr{A}}$.
$\blacksquare$
\\
In ~\cite{ZCLL}, Zhen Chen , Lin-zhang Lu and Zhi-bing Liu proposed and proved a similar property as Proposition 4.1 by a different method.
\\\par

\section{Tensor Transpose and Tensor Decomposition} \label{TD-sec}
There are a number of tensor decompositions among which CANDECOMP/PARAFAC decomposition and Tucker decomposition are most popular ~\cite{TDAA}~\cite{Comon}. They are used in psychometrics, applied statistics, weblink analysis and many other fields. In this section, we will focus on the relationship between transpose and the two major decompositions.
\subsection{CP decomposition}
CP decomposition is short for CANDECOMP/PARAFAC decomposition, which are introduced by Hitchcock~\cite{Hit1}~\cite{Hit2}, Cattell~\cite{RBC1}~\cite{RBC2}, Carroll and Chang~\cite{CC}, and Harshman~\cite{Harshman}. The CP decomposition is strongly linked with rank-one tensors. Usually, an $n$-order rank-one tensor can be written in the outer product of $n$ vectors. For example,
$\boldsymbol{\mathscr{X}}{\in}{\mathbb{R}}^{I_{1}{\times}{I_{2}}{\times}{I_{3}}}$, $\boldsymbol{\mathscr{X}}=\textbf{u}_{1}{\circ}\textbf{u}_{2}{\circ}\textbf{u}_{3}$, then $\boldsymbol{\mathscr{X}}$ is a rank-one tensor, where $\textbf{u}_{1},\textbf{u}_{2},\textbf{u}_{3}$ are vectors and $\circ$ is outer product operator.
\begin{proposition}
If rank-one tensor $\boldsymbol{\mathscr{X}}=\textbf{u}_{1}{\circ}\textbf{u}_{2}{\circ}\cdots{\circ}\textbf{u}_{N}$,  $\boldsymbol{\mathscr{X}}^{T\sigma}=\textbf{u}_{\sigma(1)}{\circ}\textbf{u}_{\sigma(2)}{\circ}\cdots{\circ}\textbf{u}_{\sigma(N)}$.
\end{proposition}
Proof: \\Assume $\boldsymbol{\mathscr{X}}{\in}{\mathbb{R}}^{I_{1}{\times}{I_{2}}{\times}\cdots{\times}{I_{N}}}$,\\
 then \begin{equation*}\boldsymbol{\mathscr{X}}=\textbf{u}_{1}{\circ}\textbf{u}_{2}{\circ}\cdots{\circ}\textbf{u}_{N}\end{equation*}
 that is \begin{equation*}\boldsymbol{\mathscr{X}}(i_{1},i_{2},\ldots,i_{N})={\textbf{u}_{1}(i_{1})}{\textbf{u}_{2}(i_{2})}\cdots{\textbf{u}_{N}(i_{N})},\end{equation*} and
 \begin{equation*}\boldsymbol{\mathscr{X}}^{T\sigma}(i_{1},i_{2},\ldots,i_{N})=\boldsymbol{\mathscr{X}}(i_{\sigma^{-1}(1)},i_{\sigma^{-1}(2)},\ldots,i_{\sigma^{-1}(N)})=
{\textbf{u}_{1}(i_{\sigma^{-1}(1)})}{\textbf{u}_{2}({i_{\sigma^{-1}(2)}})}\cdots{\textbf{u}_{N}({i_{\sigma^{-1}(N)}})}\end{equation*}
because for all $k \in \{1,2,\cdots,N\},$ there exists $j$, s.t. $\sigma(j)=k, \sigma^{-1}(k)=j$\\
so \begin{equation*} \textbf{u}_{k}(i_{\sigma^{-1}(k)})=\textbf{u}_{\sigma(k)}(i_{j})\end{equation*}\\
therefore \begin{equation*}{\textbf{u}_{1}(i_{\sigma^{-1}(1)})}{\textbf{u}_{2}({i_{\sigma^{-1}(2)}})}\cdots{\textbf{u}_{N}({i_{\sigma^{-1}(N)}})}={\textbf{u}_{\sigma(1)}(i_{1})}{\textbf{u}_{\sigma(2)}({i_{2}})}\cdots{\textbf{u}_{\sigma(N)}({i_{N}})}\end{equation*}\\
that is \begin{equation*}\boldsymbol{\mathscr{X}}^{T\sigma}=\textbf{u}_{\sigma(1)}{\circ}\textbf{u}_{\sigma(2)}{\circ}\cdots{\circ}\textbf{u}_{\sigma(N)}.\blacksquare\end{equation*}
\\\par
The CP decomposition factorizes a tensor into a sum of rank-one tensors. Take 3-order situation as an example. Let $\boldsymbol{\mathscr{X}}$ be a 3-order tensor, and $\boldsymbol{\mathscr{X}}_{i}$ be rank-one 3-order tensors, $\boldsymbol{\mathscr{X}}_{i}={\textbf{a}_{i}}\circ{\textbf{b}_{i}}\circ{\textbf{c}_{i}}$. Then the CP decomposition of $\boldsymbol{\mathscr{X}}$ can be written as
\begin{equation*}
\boldsymbol{\mathscr{X}}=\sum_{i=1}^{R}\boldsymbol{\mathscr{X}}_{i}=\sum_{i=1}^{R}{\textbf{a}_{i}}\circ{\textbf{b}_{i}}\circ{\textbf{c}_{i}}
\end{equation*}
We denote matrix $A_{1}$, $A_{2}$, $A_{3}$, as the combination of vectors $\textbf{a}_{i}$, $\textbf{b}_{i}$, $\textbf{c}_{i}$, i.e., $A_{1}=\left(
                                                                                                                               \begin{array}{cccc}
                                                                                                                                 \textbf{a}_{1} & \textbf{a}_{2} & \cdots & \textbf{a}_{R} \\
                                                                                                                               \end{array}
                                                                                                                             \right)
$. Then CP decomposition can be expressed by
\begin{equation*}
\boldsymbol{\mathscr{X}}=\sum_{i=1}^{R}{\textbf{a}_{i}}\circ{\textbf{b}_{i}}\circ{\textbf{c}_{i}}:=[[A_{1},A_{2},A_{3}]]
\end{equation*}
According to Proposition 5.1, we have a property as follows,
\begin{proposition}
If $\boldsymbol{\mathscr{X}}=[[A_{1},A_{2},A_{3}]]$, then $\boldsymbol{\mathscr{X}}^{T\sigma}=[[A_{\sigma(1)},A_{\sigma(2)},A_{\sigma(3)}]]$.
\end{proposition}
\par
The rank of a tensor is defined as the smallest number of rank-one tensors that exactly sum up to that tensor. From previous research ~\cite{TDAA}, rank decompositions are often unique. Then we have the following property.
\begin{proposition}
The rank of a certain tensor is invariant under any transpose.
\end{proposition}
\par

\subsection{Tucker decomposition}
The Tucker decomposition was first introduced by Tucker~\cite{Tucker1}~\cite{Tucker2} and it decomposes a tensor into a core tensor multiplied by a matrix along each mode. Here, we consider the 3-order situation. Let $\boldsymbol{\mathscr{X}}$ be a 3-order tensor of size $I_{1}{\times}I_{2}{\times}I_{3}$, we have
\begin{equation*}
\boldsymbol{\mathscr{X}}=\boldsymbol{\mathscr{G}}{\times_{1}}A_{1}{\times_{2}}A_{2}{\times_{3}}A_{3}
={\sum_{j_{1}}^{J_{1}}}{\sum_{j_{2}}^{J_{2}}}{\sum_{j_{3}}^{J_{3}}}g_{j_{1}j_{2}j_{3}}\textbf{a}_{j_{1}}\circ\textbf{b}_{j_{2}}\circ\textbf{c}_{j_{3}}
:=[[\boldsymbol{\mathscr{G}};A_{1},A_{2},A_{3}]],
\end{equation*}
where $A_{i}{\in}\mathbb{R}^{I_{i}{\times}J_{i}}$ are the factor matrices. The tensor $\boldsymbol{\mathscr{G}}$ of size $J_{1}{\times}J_{2}{\times}J_{3}$ is called the core tensor of $\boldsymbol{\mathscr{X}}$.
\begin{proposition}
If $\boldsymbol{\mathscr{X}}=[[\boldsymbol{\mathscr{G}};A_{1},A_{2},A_{3}]]$,
 $\boldsymbol{\mathscr{X}}^{T\sigma}=[[\boldsymbol{\mathscr{G}}^{T\sigma};A_{\sigma(1)},A_{\sigma(2)},A_{\sigma(3)}]]$.
\end{proposition}
Proof:\\
According to Proposition 3.1, we have\\ \begin{equation*}\boldsymbol{\mathscr{X}}^{T\sigma}=(\boldsymbol{\mathscr{G}}{\times_{1}}A_{1}{\times_{2}}A_{2}{\times_{3}}A_{3})^{T\sigma}=(\boldsymbol{\mathscr{G}}{\times_{1}}A_{1}{\times_{2}}A_{2})^{T\sigma}{\times_{\sigma^{-1}(3)}}A_{3}
=\boldsymbol{\mathscr{G}}^{T\sigma}{\times_{\sigma^{-1}(1)}}A_{1}{\times_{\sigma^{-1}(2)}}A_{2}{\times_{\sigma^{-1}(3)}}A_{3}\end{equation*}\\
so
\begin{equation*} \boldsymbol{\mathscr{X}}^{T\sigma}=\boldsymbol{\mathscr{G}}^{T\sigma}{\times_{\sigma^{-1}(1)}}A_{1}{\times_{\sigma^{-1}(2)}}A_{2}{\times_{\sigma^{-1}(3)}}A_{3}\end{equation*}\\
and because \begin{equation*}\boldsymbol{\mathscr{G}}{\times_{1}}A_{1}{\times_{2}}A_{2}=\boldsymbol{\mathscr{G}}{\times_{\sigma(1)}}A_{\sigma(1)}{\times_{\sigma(2)}}A_{\sigma(2)}, \text{ seen in} ~\cite{MTC},\end{equation*}\\
we have \begin{equation*}\boldsymbol{\mathscr{G}}^{T\sigma}{\times_{\sigma^{-1}(1)}}A_{1}{\times_{\sigma^{-1}(2)}}A_{2}{\times_{\sigma^{-1}(3)}}A_{3}
=\boldsymbol{\mathscr{G}}^{T\sigma}{\times_{1}}A_{\sigma(1)}{\times_{2}}A_{\sigma(2)}{\times_{3}}A_{\sigma(3)},\end{equation*}\\
therefore \begin{equation*}\boldsymbol{\mathscr{X}}^{T\sigma}=[[\boldsymbol{\mathscr{G}}^{T\sigma};A_{\sigma(1)},A_{\sigma(2)},A_{\sigma(3)}]]. \blacksquare\end{equation*}\\\par
From this property, we can see that $\boldsymbol{\mathscr{G}}^{T\sigma}$ is the core tensor of $\boldsymbol{\mathscr{X}}^{T\sigma}$, if $\boldsymbol{\mathscr{G}}$ is the core tensor of $\boldsymbol{\mathscr{X}}$. \par
Actually, Proposition 5.2, 5.3, 5.4 also apply to 4 or higher order tensors.

\section{Conclusion}
The notion of tensor transpose is often mentioned together with supersymmetric tensors but specific discussions concerning tensor transpose draw less attention. In this paper, we propose the definition of tensor transpose and proved some basic properties. According to Proposition 3.1, properties regarding inner product, tensor eigenvalues, tensor decompositions and rank are derivated in following sections. In future work, the introduction of tensor transpose may be useful for tensor theory research, computation or algorithms improvement. We will keep working on it.

\end{document}